\documentclass[a4paper,12pt]{article}
\usepackage[english]{babel}
\usepackage{amsmath,amssymb,amsthm}
\usepackage{hyperref}
\usepackage{comment}

\date{}

\newtheorem{theorem}{Theorem}

\newtheorem{proposition}{Proposition}
\newtheorem{corollary}{Corollary}
\theoremstyle{definition}
\newtheorem{problem}{Problem}
\newtheorem{remark}{Remark}

\title{Projective tilings and full-rank perfect codes%
\thanks{The work is funded by the Russian Science Foundation, Grant 22-11-00266.
}
\thanks{This version of the article has been accepted for publication, after peer
review but is not the Version of Record and does not reflect post-acceptance
improvements. The Version of Record is available online at:
\url{https://doi.org/10.1007/s10623-023-01256-y}. Use of this Accepted Version is subject to the publisher's Accepted
Manuscript terms of use \url{https://www.springernature.com/gp/open-research/policies/accepted-manuscript-terms}. 
}
}
\author{%
 Denis S. Krotov%
\thanks{Sobolev Institute of Mathematics, Novosibirsk 630090, Russia}
}

 \newcommand\plangle{\langle\hspace{-0.6ex}\langle}
 \newcommand\prangle{\rangle\hspace{-0.6ex}\rangle}
\newcommand\alangle{\langle\hspace{-0.5ex}\langle}
\newcommand\arangle{\rangle\hspace{-0.5ex}\rangle}
\def\vc#1{\bar{#1}}
\def\pp#1{\dot{#1}}
\def\ppp#1{\ddot{#1}}
\def\FF{\mathbb{F}}
\def\ZZ{\mathbb{Z}}
\def\BB{\mathcal{B}}
\def\rank{\mathrm{rank}}
\def\ker{\mathrm{ker}}
\def\per{\mathrm{per}}
\def\wt{\mathrm{wt}}

\def\PG{\mathrm{PG}}
\def\SS{\mathbb{S}}

\begin{document}
\maketitle

\begin{abstract} A tiling of a vector space $S$ is the pair $(U,V)$ of its subsets such that every vector in $S$ is uniquely represented as the sum of a vector from $U$ and a vector from $V$. A tiling is connected to a perfect codes if one of the sets, say $U$, is projective, i.e., the union of one-dimensional subspaces of $S$. A tiling $(U,V)$ is full-rank if the affine span of each of $U$, $V$ is $S$. For finite non-binary vector spaces of dimension at least $6$ (at least $10$), we construct full-rank tilings $(U,V)$ with projective $U$ (both $U$ and $V$, respectively). In particular, that construction gives a full-rank ternary $1$-perfect code of length $13$, solving a known problem. We also discuss the treatment of tilings with projective components as factorizations of projective spaces.

Keywords: perfect codes, tilings, group factorization, full-rank tilings, projective geometry
\end{abstract}

\section{Introduction}\label{s:intro}
We are interested in factorizations 
(also known as tilings)
of an elementary
abelian group $\ZZ_p^n$ into the direct sum of two subsets, called tiles,
and the connection of such factorizations with $1$-perfect codes in Hamming spaces. 
This connection has been extensively studies
in the binary case, where each tiling~$(U,V)$
of $\ZZ_2^n$ corresponds 
to a $1$-perfect binary code of length $|U|-1$.
To establish a similar correspondence 
(see Proposition~\ref{p:1tiling} below)
for an arbitrary prime~$p$,
one should require an additional
property of $U$ to be a $\ZZ$-set,
in the sense of~\cite{Okuda}: 
$\ZZ_p\cdot U = U$.
We study factorizations $(U,V)$ of $\ZZ_p^n$ 
with $U$ or both $U$ and $V$ 
satisfying this property or,
more generally,
factorizations $(U,V)$ of $\FF_q^n$ (where 
$\FF_q$ is the finite field 
of a prime power order $q$) 
such that $\FF_q\cdot U = U$.

A notable particular consequence of the results of this paper
is proving the existence of a ternary $1$-perfect code of 
length~$13$ and rank~$13$, which was the last open question 
in the story of the study of possible values of the rank of 
$1$-perfect codes. The next paragraph contains a brief survey
of results in this area.

One-perfect codes in $\FF_q^n$ ($q$ prime power)
exist if and only if $n = (q^m-1)/(q-1)$, $m\ge 1$.
The case $m=1$ is trivial and usually not considered.
In general, the rank of a perfect code can possess
any value from $n-m$ (linear and affine codes) to $n$ (full-rank codes);
the exceptions are $1$-perfect binary length-$3$ and length-$7$ codes
and ternary length-$4$ codes, which are all of rank $n-m$.
% \begin{itemize}
%  \item 
 The construction of $1$-perfect linear codes 
 (i.e., of minimum possible rank $n-m$)
 originates to the work of Hamming~\cite{Hamming:50},
 while in the most general form, for every $q$ and $m$, 
 it was suggested by Golay~\cite{Golay:49}.
 For $q=2$, $m\le 3$~\cite{Zaremba:52} and $q=3$, $m=2$,
 $1$-perfect codes are unique up to equivalence and affine,
 which means in our context that there are 
 no codes of non-minimum rank for those parameters.
%  \item 
 The first known construction of non-linear 
 (i.e., of rank larger than the minimum possible value)
 binary $1$-perfect codes, for any $m\ge 4$,
 was proposed by Vasil'ev~\cite{Vas:nongroup_perfect.de};
 Sch\"onheim~\cite{Schonheim68} 
 generalized Vasil'ev's construction
 to any $q\ge 3$ and $m\ge 3$.
%  \item 
 In~\cite{EV:94},  Etzion and Vardy 
 constructed binary $1$-perfect
 codes of all admissible ranks, $m\ge 4$.
%  \item
 Fraser and Gordon showed in~\cite{FraGor:79} that
 permuting the symbols of the alphabet
 in one position of all codewords
 keeps the parameters of the code
 but can change its rank if $q\ge 4$
 ($4$ is the first value of~$q$ for which
 there are non-affine permutations of~$\FF_q$).
 This allows to construct $1$-perfect codes
 of any admissible rank if $q\ge 4$ and $m\ge 2$;
 a code constructed in this way is equivalent
 to a linear code if we naturally consider
 the alphabet permutations as equivalence operations.
%  \item 
 In~\cite{PheVil:2002:q},
 Phelps and Villanueva constructed
 $1$-perfect codes
 of all admissible ranks   
 for $q\ge 3$ and $m\ge 4$~\cite{PheVil:2002:q} 
 and for non-prime $q\ge 4$ and $m\ge 3$~\cite{PheVil:q-kernel}. 
%  \item 
 After the results above, only the existence of ternary 
 $1$-perfect codes of length~$\frac{3^3-1}{2-1}=13$ and 
 ranks~$12$ and~$13$ remained open.
%  (the Hamming code has rank~$10$, and a rank~$11$ code 
%  is given by the construction in~\cite{Schonheim68}).
 The rank~$12$ case was solved by Romanov~\cite{Romanov:2019}.
 The current study presents, among other results,
 a solution for the last open case.
% \end{itemize}
 \begin{remark} 
The switching technique used in 
\cite{Vas:nongroup_perfect.de},
\cite{Schonheim68},
\cite{EV:94},
\cite{PheVil:2002:q},
\cite{PheVil:q-kernel}
 results in codes nonequivalent to a linear code in general.
 However, it remains unanswered 
 whether we can reduce the rank of a code constructed
 there,
 as well as in the current paper if $q\ge 4$,
 by equivalence operations 
 (probably, not);
 this question requires further studying.
 The following fact might be important in this context:
 if a code~$D$ is obtained from a linear code~$C$ by switching
 and the number of switched codewords
 is relatively small, then the original linear code~$C$
 can be uniquely reconstructed from~$D$
 as the closest code that is equivalent to a linear code~\cite{AvgGor:2019}.
 It is also notable 
 that no $1$-perfect codes 
 that are inequivalent to linear codes
 are known if $q$ is prime and $m=2$
 (for $q\le 7$, such codes do not exist~\cite{KKO:smallMDS}).
 \end{remark}
 
In this paper, after Section~\ref{s:def}
with preliminaries, we obtain the following
results.

Firstly, 
in Section~\ref{s:tiling}, for any field $\FF$ of order larger than $2$, 
we construct full-rank tilings
$(U,V)$ of $\FF^n$
such that $\FF\cdot U = U$ ($n\ge 6$ even)
or both $\FF\cdot U = U$
and  $\FF\cdot V = V$ hold ($n\ge 10$ even,
which improves a construction for $n\ge 12$ in
\cite{Szabo:2013}).

Secondly, 
in Section~\ref{s:codes},
we construct full-rank $1$-perfect codes 
over $\FF_q$, $q\ge 3$,
in particular, 
a full-rank ternary $1$-perfect code of length~$13$.
Extending that theoretical result with
a computer-aided search, we find such 
codes with different dimension of the kernel
(the set of all periods of the code; constructing codes with different pairs [rank, kernel dimension] is also a known direction in the study of perfect codes~\cite{EV:98}, \cite{ASH:RankKernel}, \cite{PheVil:RankKernel}).
% It should be noted that,
% as was observed in~\cite{FraGor:79},
% for larger $q$ it is relatively easy to construct
% full-rank $1$-perfect codes using a simple operation, permutation of the alphabet,
% that is neither linear nor affine for $q\ge 4$,
% but preserves the Hamming distance 
% (and hence the parameters of a code).
% However, for $q\le 3$, every permutation
% of $\FF_q$ is affine, and that is why
% the last parameters for which
% the problem of existence of full-rank
% $1$-perfect codes remained unknown were
% those of ternary $1$-pefect codes of length~$13$,
% solved in the present research.

Additionally, in Section~\ref{s:factor}
we discuss the connection 
of projective tilings with 
factorizations
of projective geometries.

For more reading about factorizations
of abelian groups, see~\cite{SzaSan:factoring};
for connections of tilings of binary
spaces with perfect codes, 
see~\cite{CLVZ:tiling,EV:98};
for the problem of existence of full-rank
factorizations of abelian groups, see~\cite{Dinitz:06}.

\section{Preliminaries}\label{s:def}

Let $\SS$ be a vector space over a field $\FF$. A pair $(U,V)$ of subsets of~$\SS$
is called a \emph{tiling} (of $\SS$)
if every element of $\SS$ is uniquely represented 
as the sum $\vc{u}+\vc{v}$, $\vc{u}\in U$, $\vc{v}\in V$.
The components $U$ and $V$ of a tiling are called \emph{tiles}.
We are mostly aimed of constructing 
tilings and codes in the vector space $\FF_q^n$ 
of the $n$-tuples over the finite field $\FF_q$ of order~$q$.
However, some of the definitions and results have sense 
for infinite vector spaces as well 
(see, e.g.,~\cite{Malyugin:2019}, \cite{Malyugin:2020} for examples of infinite perfect codes).
If the reader is not interested in abstract infinite generalizations,
it is recommended to think that $\FF$ is always finite.

The \emph{rank}, $\rank(U)$, 
of a subset $U$ of a vector space $\SS$ is
the dimension of the affine span $\alangle U \arangle$ of~$U$ 
(we use double and single brackets to distinguish between the affine span
$ \alangle \,\cdot\, \arangle$ and the linear span $ \langle \,\cdot\, \rangle$).
If $\alangle U \arangle=\SS$, then $U$ 
is said to be \emph{full-rank}.
A tiling $(U,V)$ is called \emph{full-rank}
if both~$U$ and~$V$ are full-rank.

A \emph{period} of a subset~$U$ of~$\SS$ 
is an element~$\vc{v}$ of~$\SS$ 
such that $U+\vc{v} = U$.
The set of all periods of~$U$
is denoted by~$\per(U)$.
(If $\FF = \FF_q$ and $p$ is the prime divisor of~$q$,
then $\per(U)$ is necessarily 
an $\FF_p$-subspace of~$\FF_q^n$ and 
is called the \emph{$p$-kernel}
of~$U$.)
% In this case, by the \emph{dimension} of $\dim(\per(U))$ of $\per(U)$,
% we will mean the value $\frac{\dim_p(\per(U))}{\log_pq}$,
% where $\dim_p$ is the dimension over $\FF_p$.
The maximal subspace in the set of all periods of~$U$ 
is called its \emph{kernel},
$\ker(U)$ (also known as the $q$-\emph{kernel} if $\FF = \FF_q$).

A set $U$ is called \emph{aperiodic} 
if $\per(U)=\{\vc0\}$, where $\vc0$
the all-zero vector. 
A tiling $(U,V)$ is called \emph{aperiodic}
if both $U$ and $V$ are aperiodic.

A subset $U$ 
of $\SS$ is called \emph{projective}
if it is the union of $1$-dimensional subspaces of $\SS$,
i.e., if $\FF \cdot U = U$.
A tiling $(U,V)$ is called \emph{semiprojective}
if $U$ is projective.
A tiling $(U,V)$ is called \emph{projective}
if both $U$ and $V$ are projective.

In the definitions above, $\SS$ can be considered as an abstract vector space 
and the basis does not play any role. 
This is not the case in the following group of definitions, 
where we consider 
the vector space $\FF^n$ of $n$-tuples over $\FF$,
considered as vectors written is some fixed basis.

The Hamming \emph{weight} $\wt(\vc{x})$ 
of a vector $\vc{x}$ in $\FF^n$
is the number of its nonzero coordinates.
The Hamming \emph{distance} between two vectors 
$\vc{x}$, $\vc{y}$ in $\FF^n$
is the number of coordinates in which they differ, i.e., 
$\wt(\vc{y}-\vc{x})$.
The radius-$r$ (Hamming) \emph{ball}
$\BB_r(\vc{x})$
centered in $\vc{x}\in\FF^n$
is the set of vectors at distance 
at most~$r$ from $\vc{x}$; we denote 
$\BB_r=\BB_r(\vc0)$.

A subset $C$ of $\FF^n$ is called an \emph{$r$-perfect code}
if $(\BB_r,C)$ is a tiling.

Nontrivial ($0<r<n$) perfect codes
over finite fields $\FF_q$ are known to exist if and only if $q=2$,
$r=\frac{n-1}2$, odd $n$ (binary repetition codes),
 $q=2$, $r=3$, $n=23$ (the binary Golay code),
 $q=3$, $r=2$, $n=11$ (the ternary Golay code),
 and $r=1$, $n=\frac{q^m-1}{q-1}$, $m=2,3,\ldots$, 
 for any prime power~$q$ \cite{Tiet:1973,ZL:1973}.
 So, all parameters of perfect codes (over fields)
 are known; however, there are many nonequivalent $1$-perfect codes \cite{Vas:nongroup_perfect.de,Schonheim68} and the class
 of all such codes is hardly possible to characterize
 in any constructive terms.
 
 The following proposition connects tilings and perfect codes.
Its first claim generalizes~\cite[Propositions~7.1]{CLVZ:tiling}, 
where it was shown for $\FF_2$ (when all tilings are projective),
 and also is a special case of~\cite[Theorem~2.1]{BloLam:84},
 where more general objects, coverings, are considered.
 
 \begin{proposition}\label{p:1tiling}
  Let $(U,V)$ be a semiprojective tiling of a vector space $\SS$ 
  over a field~$\FF$.
  Let $U^*$ be a complete set
  of 
  mutually non-colinear
  representatives of $U$.
  We assume that $U^*$ is finite and denote $N=|U^*|$
  {\rm(}in particular, $N\cdot(|\FF|-1)=|U|-1$${\rm)}$.
  Let the matrix~$H$
  be formed by the elements of~$U^*$ as columns.
  If 
  \begin{equation}\label{eq:C}
  C = \{ c \in \FF^{N}:\ Hc  \in V \}, 
  \end{equation}
  then
  \begin{itemize}
   \item[\rm(i)] $C$ is a $1$-perfect code in $\FF^N$,
   \item[\rm(ii)]  $\rank(C) = \rank(V_U) + N -  r $,
   \item[\rm(iii)]  $\dim(\ker(C)) = \dim(\ker(V_U)) + N -  r $,\quad
     $|\per(C)| = |\per(V_U)|\cdot |\FF|^{N - r} $,
  \end{itemize}
  where $V_U = V\cap \langle U \rangle$
and $ r  = \rank(U)$
  {\rm(}if $U$ is full-rank, then
   $V_U = V$ and 
  $ r =\dim(\SS)${\rm)}.
 \end{proposition}
 
\begin{proof}
\newcommand\UU{\overline{U}}
Denote $\UU = \langle U \rangle$.
Since $(U,V)$ is a tiling,
every element of $\UU$ is uniquely represented 
as the sum of some $\vc{u}$ from $U$ and $\vc{v}$ from $V$;
moreover, trivially $\vc{u}\in \UU$ 
and hence $\vc{v}\in \UU$
and  $\vc{v}\in V_U$.
By the definition, $(U,V_U)$ is a tiling of $\UU$.

Moreover, \eqref{eq:C} is equivalent to
$C = \{ c \in \FF^{N}:\ Hc  \in V_U \}$ 
because $Hc$ belongs to $\UU$
and cannot belong to $V \setminus \UU$.

To prove (i), we consider an arbitrary $\vc{x}$ from $\FF^N$
and denote $\vc{s} = H\vc{x}$. 
Since $\vc{s} \in \UU$ and 
$(U,V_U)$ is a tiling of $\UU$,
we find that $\vc{s}$ is 
uniquely 
represented
in the form $\vc{s} = \vc{u} + \vc{v}$, 
$\vc{u}\in U$,  $\vc{v}\in V_U$.
If $\vc{u} = \vc{0}$, then $\vc{u} = H \vc{0} $;
otherwise, by the definition of the matrix~$H$,
$\vc{u}$ is uniquely represented as
a multiple of one of its columns,
i.e.,  $\vc{u} = H \vc{b} $ for some~$\vc{b}$ 
in $\BB_1\backslash\{\vc{0}\}$.
In both cases, $\vc{u} = H \vc{b} $ for a unique~$\vc{b}$ 
in $\BB_1$.
Denoting $\vc{c} = \vc{x} - \vc{b}$,
we see that 
$ H\vc{c} = \vc{y}-\vc{u} = \vc{v} \in V_U$
and hence $\vc{c} \in C$.
So,
$\vc{x} = \vc{c} + \vc{b}$, where $\vc{c} \in C$ and $\vc{b}\in \BB_1$.
We also see $\vc{b}$ from $\BB_1$  and hence $\vc{c}$ from $C$
are uniquely determined from $\vc{x}$. 
Therefore, by the definitions,
$(\BB_1,C)$ is a tiling and $C$ is a $1$-perfect code.

Next, $H$ acts as a linear operator from 
$\FF^N$ onto $\UU$.
Hence, it maps $\alangle C \arangle$ to $\alangle V_U \arangle$.
Since the preimage (under $H$) of every vector in $\UU$
 is an affine subspace of dimension ${N -  r }$, 
we find 
$\dim(\alangle C \arangle) = 
 \dim(\alangle V_U \arangle) + ( N - r )$
and $\rank(C) = \rank(V_U) + N -  r $, which proves~(ii).

(iii) is
% (iii) and (iv) are 
straightforward from the obvious fact that
multiplication by~$H$ maps 
$\ker(C)$ to $\ker(V_U)$ 
and
$\per(C)$ to $\per(V_U)$.
\end{proof}

\begin{remark}
 The theorem above shows the connection of $1$-perfect codes 
and projective tilings. 
Tilings of a vector space over
$\FF=\FF_q$, where $q$ is a power of a prime~$p$,
are essentially factorizations of an elementary $p$-group.
As shown in~\cite{Wu:additive}, some perfect codes can be constructed
from factorizations of non-elementary abelian $p$-groups.
\end{remark}

\section{Constructing full-rank
semiprojective and projective tilings}\label{s:tiling}

The construction of a tiling in the proofs 
of the following two theorems
exploits the idea from Step~1 of the proof 
of~\cite[Theorem~2]{Szabo:2006}
(see also \cite[Theorem~2.5]{denBreeijen:master}, 
where that proof is given in a notation close to ours 
and for any $m\ge 3$), but we modified it to make $U$ 
(both $U$ and $V$, in Theorem~\ref{th:n8})
projective.

\begin{theorem}\label{th:n6}
 If $\FF$ is a field of cardinality $|\FF|$ larger than $2$
 and $m$ is an integer, $m\ge 3$, 
 then there is a full-rank aperiodic semiprojective tiling $(U,V)$ of 
 $\FF^{2 m}$ with $|U|=|V|=|\FF|^{m}$.
\end{theorem}

\begin{proof}
 Let 
 $( \vc{x}_1, \ldots , \vc{x}_m,
    \vc{y}_1, \ldots , \vc{y}_m )$
 be the natural basis of $\FF^{2 m}$. 
 For convenience, we identify 
 $\vc{x}_{m+i}=\vc{x}_{i}$ and 
 $\vc{y}_{m+i}=\vc{y}_{i}$.
 We start with constructing a non-full-rank semiprojective tiling
 $(H,V)$.
 Let $H = \langle \vc{x}_1, \ldots , \vc{x}_m \rangle$,
 and let $V = V _1 + \ldots + V _m$, where
 $$
 V _i = 
%  \{ \gamma \vc{y}_i: \ \gamma \in \FF\backslash\{1\} \}
%  \cup \{\vc{y}_i+\vc{x}_i\} 
\big( \langle \vc{y}_i \rangle 
 \setminus \{\vc{y}_i\} \big)
 \cup  \{\vc{y}_i+\vc{x}_i\}
 .
 $$
 It is not difficult to see that $(H,V)$ is a tiling 
 (indeed, $V$ is obtained from $\langle \vc{y}_1, \ldots , \vc{y}_m \rangle$
 by adding periods of $H$ to some elements).
 It is easy to see that $V$ is full-rank and aperiodic
 (see (iii) and (v) below), and
 it remains to modify $H$ 
 with making it full-rank and aperiodic too (but keeping projective).
 The modification is based on the following fact.
 
 (*) \emph{If $L$ is a subset of $\FF^{2 m}$ with a period $\vc{x}_i$,
 then $L+V = (L+\gamma\vc{y}_i) + V$ for every $\gamma$ in $\FF$.}
 To prove it, we assume $i=1$ 
 without loss of generality.
 Since $\vc{x}_1$ is a period of $L$, subtracting it from some elements
 of $V$ does not change the sum $L+V$. Hence, 
 \begin{equation}\label{e:LVLV'}
  L+V=L+V',
 \end{equation}
 where $V' = \langle\vc{y}_1\rangle + V_2 + \ldots + V _m$.
 Now we see that $\gamma \vc{y}_1$ is a period of $V'$,
 and so 
 \begin{equation}\label{e:LVLyV'}
 L+V'= L+\gamma\vc{y}_1+V'.
 \end{equation}
 Again, $\vc{x}_1$ is a period of $L+\gamma\vc{y}_1$,
 and 
 \begin{equation}\label{e:LyV'LyV}
 L+\gamma\vc{y}_1+V'= L+\gamma\vc{y}_1+V, 
 \end{equation}
  Equalities~\eqref{e:LVLV'}--\eqref{e:LyV'LyV} prove~(*).
 
 We construct $U$ by modifying $H$ in the following manner. 
 For $i \in \{1, \ldots ,m\}$ and $\gamma \in \FF \backslash \{0\}$,   
 denote 
 \begin{eqnarray*}
   H_{i,\gamma} &=& 
 \langle \vc{x}_i\rangle + \gamma \vc{x}_{i+1}
 \\\
 \mbox{and}
 \qquad
 U_{i,\gamma} &=& 
 \langle \vc{x}_i\rangle + \gamma \vc{y}_{i} + \gamma \vc{x}_{i+1};
 \end{eqnarray*}
 it is not difficult to observe that all these sets are mutually disjoint.
 Then, define
 $$
 U = 
 H 
 \setminus 
 \Big(\bigcup_{i,\gamma} H_{i,\gamma}\Big)
 \cup
 \Big(\bigcup_{i,\gamma} U_{i,\gamma}\Big) 
 $$
 (the unions are over $i \in \{1, \ldots ,m\}$ and $\gamma \in \FF \backslash \{0\}$).
 Now, we can claim that $(U,V)$ is a required full-rank aperiodic semiprojective
 tiling.

 (i) \emph{$U$ is projective} because $H$, 
 $\{\vc0\}\cup \bigcup_{\gamma\in\FF\backslash\{0\}} H_{i,\gamma}$,
 and
 $\{\vc0\}\cup \bigcup_{\gamma\in\FF\backslash\{0\}} U_{i,\gamma}$
 are projective, $i=1, \ldots ,m$.
 
 (ii) \emph{$U$ is full-rank}. Indeed, it is easy to see
 that 
 $
 H 
 \setminus 
 \big(\bigcup_{i,\gamma} H_{i,\gamma}\big)
 $
 still spans $H$ 
 (it includes the subset 
 $(\FF \backslash \{0\})^m$, for example), while each $U_{i,1}$
 adds $\vc{y}_i$ to the span, $i=1, \ldots ,m$.
 
 (iii) \emph{$V$ is full-rank} because 
 $\alangle V_i \arangle=\langle \vc{x}_i,\vc{y}_i \rangle$, $i=1, \ldots ,m$.
 
 (iv) \emph{$U$ is aperiodic}. Indeed, consider the intersections
 of $U$ with the cosets of the subspace $H$.
 There are three types of such intersections:
 $ U \cap H $;
 $ U \cap (H + \gamma \vc{y}_i) = U_{i,\gamma}$, 
 $i \in \{1, \ldots ,m\}$, $\gamma \in \FF \backslash \{0\}$;
 and the empty intersections.
 Since $ U \cap H $ is the only intersection that is neither a line nor empty,
 we see that $U$ has no periods out of $H$.
 On the other hand, the two lines 
 $U_{1,1} = U\cap (H+\vc{y}_1)$ and $U_{2,1} = U\cap (H+\vc{y}_2)$ 
 has no common nonzero periods; 
 hence, $U$ has no nontrivial periods in $H$.
 
 (v) \emph{$V$ is aperiodic}. Indeed, the coset 
 $\vc{x}_1 + \ldots +\vc{x}_m + Y$ of
 $Y = \langle \vc{y}_1 , \ldots ,\vc{y}_m \rangle$ contains only one 
 element of $V$; hence, there are no nontrivial periods of $V$ in $Y$.
 On the other hand, there are no other cosets of $Y$ that intersect with $V$
 in exactly one element; this certifies that $V$ has no periods out of $Y$.
 
 (vi) \emph{$(U,V)$ is a tiling}. Indeed, $(H,V)$ is a tiling,
 and from (*) we see $H_{i,\gamma}+V = U_{i,\gamma}+V$. 
 Hence, replacing $H_{i,\gamma}$ by $U_{i,\gamma}$ does not change the
 tiling property.
\end{proof}

It is quite natural to ask if the construction 
can be strengthened with making both tiles projective.
The following variant of the construction do
this for $m\ge 5$, while the existence of full-rank 
projective tilings of $\FF_q^n$ for $n<10$
remains an open challenging problem.
Another construction of full-rank 
projective tilings of $\FF_q^{2m}$ for $m\ge 6$
was suggested in~\cite{Szabo:2013}.

\begin{theorem}\label{th:n8}
 If $\FF$ is a field of cardinality $|\FF|$ larger than $2$
 and $m$ is an integer, $m\ge 5$, 
 then there is a full-rank aperiodic projective tiling $(U,V)$ of 
 $\FF^{2 m}$ with $|U|=|V|=|\FF|^{m}$
\end{theorem}

\begin{proof}
The idea is the same as 
in the proof of Theorem~\ref{th:n6}, 
but now
we cannot represent $V$ as the sum 
$V_1 + \ldots + V_m$ 
(indeed, if all $V_i$ are projective, then such
sum is a proper 
subspace of $\FF^{2m}$ and $V$ is not full-rank).
Let $B = \langle \vc{y}_1, \ldots , \vc{y}_m \rangle$,
and we define 
$$ V = \{ v(\vc{z}):\ \vc{z} \in B \}, $$ 
where the function $v(\cdot)$ is defined as
$$
v(0, \ldots , 0, z_1, \ldots , z_m ) = 
((z_2{=}0)?z_1{:}0,\  (z_3{=}0)?z_2{:}0,\ \ldots ,\ (z_1{=}0)?z_m{:}0,\ z_1, \ldots , z_m),
$$
where $\nu?\lambda{:}\mu$ is a short notation
for ``if $\nu$ then $\lambda$ else $\mu$''.
% (Remark: if we replace the conditions $(z_2{=}0)$, $(z_3{=}0)$, \ldots
%  by $(z_1{=}1)$, $(z_2{=}1)$ \ldots, respectively, then we obtain $V$
%  from the proof of Theorem~\ref{th:n6}.)
It is not difficult to see that for a $1$-dimensional subspace $L$ of $B$,
its image $v(L)$ is also a $1$-dimensional space
(indeed, for each $i$, the condition $z_i=0$ is either constantly true 
or constantly false on all non-zero elements of $L$, and so $v(\cdot)$
 is just a linear function on $L$). Hence, $V$ is projective.
 
 To define $U$, we keep the notation $H=\langle\vc{x}_1, \ldots,  \vc{x}_m\rangle $ and
 $$
 U = 
 H 
 \setminus 
 \Big(\bigcup_{i,\gamma} H_{i,\gamma}\Big)
 \cup
 \Big(\bigcup_{i,\gamma} U_{i,\gamma}\Big) 
%  \quad \mbox{(the unions are over $i \in \{1, \ldots ,m\}$ and $\gamma \in \FF \backslash \{0\}$)}
,
 $$
 but now
\begin{eqnarray*}
 H_{i,\gamma}&=&\langle \vc{x}_i, \vc{x}_{i+1} \rangle
                  + \gamma \vc{x}_{i+2}    
\\
\mbox{and} 
\quad
U_{i,\gamma}&=&\langle \vc{x}_i, \vc{x}_{i+1} \rangle
                  + \gamma \vc{x}_{i+2}
                  + \gamma \vc{y}_{i+1}.
\end{eqnarray*}
For $m\ge 5$, all  $H_{i,\gamma}$ and $U_{i,\gamma}$
are mutually disjoint 
(which is not true for $m=3,4$, 
and this is the only reason 
why we require $m\ge 5$).
The key observation to prove that $(U,V)$
is a tiling is now the following:

 (*) \emph{if $L$ is a subset of $\FF^{2m}$ with  periods $\vc{x}_i$ and $\vc{x}_{i+1}$,
 then $L+V = (L+\gamma\vc{y}_{i+1}) + V$ for every $\gamma$ in $\FF$.} 
  The proof of (*) essentially repeats the similar 
  for Theorem~\ref{th:n6}.
  We assume $i=1$.
 Since $\vc{x}_1$ and $\vc{x}_2$ are periods of $L$,
 varying the values in the first two coordinates
 of vectors of $V$ 
 does not change the sum $L+V$. Hence, 
 \begin{equation}\label{e:LVLV'_}
  L+V=L+V',
 \end{equation}
 where $ V' = \{ v'(\vc{z}):\ \vc{z} \in B \}, $ 
$$
v'(0, \ldots , 0, z_1, \ldots , z_m ) = 
(0,\,0,\  (z_4{=}0)?z_3{:}0,\ \ldots ,\ (z_1{=}0)?z_m{:}0,\ z_1,z_2,z_3, \ldots , z_m).
$$
 Now, varying the value of $z_2$, 
 we see that $\gamma\vc{y}_2$ is a period of $V'$,
 and so 
 \begin{equation}\label{e:LVLyV'_}
 L+V'= L+\gamma\vc{y}_2+V'.
 \end{equation}
 Again, $\vc{x}_1$ and  $\vc{x}_2$ 
 are periods of $L+\gamma\vc{y}_2$,
 and 
 \begin{equation}\label{e:LyV'LyV_}
 L+\gamma\vc{y}_2+V'= L+\gamma\vc{y}_2+V.
 \end{equation}
 Equalities~\eqref{e:LVLV'_}--\eqref{e:LyV'LyV_} prove~(*).
 
 The rest of the proof of the theorem
is similar to that of Theorem~\ref{th:n6}.
\end{proof}

\section{Full-rank perfect codes}\label{s:codes}

From Theorem~\ref{th:n6} and Proposition~\ref{p:1tiling}, we have the following.

\begin{corollary}\label{c:r13}
 In $\FF_q^n$, $n=\frac{q^m-1}{q-1}$, $q\ge 3$, $m\ge 3$, 
 there is a full-rank $1$-perfect 
 code with the dimension $n-2m$ of the kernel. 
\end{corollary}

As was noted in the introduction,
the existence of full-rank $1$-perfect codes in 
$\FF_3^{13}$ was an open problem, and by Corollary~\ref{c:r13}
we can construct such code with 
 kernel dimension~$7$.
Starting with that code and 
using the switching approach,
full-rank $1$-perfect codes in $\FF_3^{13}$
with 
 kernel dimensions~$6$, $5$, $4$, and~$3$ were found
 by computer search
 (the two-coordinate switchings described in~\cite{KrotovShi:3ary} was applied,
 while the traditional one-coordinate switching did not work
 for the given starting code in the meaning that the switched code was always
 equivalent to the original one). 
 
\begin{theorem}\label{th:r13}
 For each $k$ from $3$, $4$, $5$, $6$, $7$,
 there is a full-rank ternary $1$-perfect code of length $13$ with 
 kernel dimension $k$.
\end{theorem}
The examples of codes are available in~\cite{Perfect-related}. 
% \url{https://dx.doi.org/10.21227/w856-4b70}.
 As follows from~\cite{OstSza:2007}, $7$ is the maximum value of kernel dimension
 for full-rank $1$-perfect codes in $\FF_3^{13}$,
 and so only the existence of codes of kernel dimension less than $3$
 remains unsolved for these parameters 
 (both for full-rank codes and for codes of rank $12$, see the rank--kernel table in~\cite{KrotovShi:3ary}).

\section{Factorizations in point-line geometries}\label{s:factor}

\newcommand{\GG}{{\mathbb G}}

This section does not contain any results apart of
presenting a geometrical treatment of projective tilings 
in terms of factorizations of projective spaces.
This interesting reformulation seems to be never observed before.

Let $\GG$ be 
% an arbitrary point--line incidence structure (geometry), for example, 
a projective space.
We will say that
two disjoint sets
$\mathcal{U}$, $\mathcal{V}$
of points of $\GG$ form a \emph{factorization}
$(\mathcal{U},\mathcal{V})$
of $\GG$ if for every point $\pp{x}$ not in 
$\mathcal{U} \cup \mathcal{V}$ there are unique
points 
$\pp{u}$ in $\mathcal{U}$
and
$\pp{v}$ in $\mathcal{V}$ such that
$\pp{x} \in \plangle \pp{u},\pp{v}\prangle$,
where $\plangle\, \cdot\, \prangle$
is the minimum projective subspace
that contains
the points inside the brackets
(in particular, $\plangle \pp{u},\pp{v}\prangle$ is a line if 
$ \pp{u} \neq \pp{v}$).

A set $\mathcal{U}$ of points is \emph{full-rank}
if $\plangle \mathcal{U} \prangle = \GG$.
A factorization $(\mathcal{U}$, $\mathcal{V})$
is \emph{full-rank} if both 
$\mathcal{U}$ and $\mathcal{V}$ are full-rank.

A point $\pp{u}$ is called a \emph{period}
of a set~$\mathcal{U}$ if  $\mathcal{U}$
is the union of lines through~$\pp{u}$.
A point set is \emph{aperiodic}
if it has no periods. 
A factorization $(\mathcal{U}$, $\mathcal{V})$
is \emph{aperiodic}
if both $\mathcal{U}$ and $\mathcal{V}$
 are aperiodic.

The following proposition connects
factorizations of projective spaces with projective tilings of vector spaces.
For a vector space $\SS$, by $\PG(\SS)$ we denote the projective space 
(of dimension $\dim(\SS)-1$)
whose points are the $1$-dimensional subspaces of $\SS$
and lines corresponding to $2$-dimensional subspaces of $\SS$.

\begin{proposition}
 Let $U$ and $V$ be two sets of points of
a vector space $\SS$. Assume that 
 $U$ and $V$ are projective, i.e.,
 $U =\bigcup_{\pp{u}\in\mathcal{U}}\pp{u}$
 and 
 $V =\bigcup_{\pp{v}\in\mathcal{V}}\pp{v}$
 for some sets $\mathcal{U}$ and $\mathcal{V}$
 of $1$-dimensional subspaces of $\SS$.
 \begin{itemize}
 
  \item[\rm(i)] $(U,V)$ is a tiling of  $\SS$
 if and only if $(\mathcal{U},\mathcal{V})$
 is a factorization of $\PG(\SS)$.
 \item[\rm(ii)]
 $U$ {\rm(}similarly, $V${\rm)} is full-rank if and only
 if $\mathcal{U}$ is full-rank.
 \item[\rm(iii)] There is a period in $\mathcal{U}$ 
 {\rm(}similarly, in $\mathcal{V}${\rm)} if and only if
 $\dim(\ker(U))\ge 1$ {\rm(}in the case of a prime-order field, it is sufficient to say that 
 $U$ has at least one nonzero period{\rm)}.
 \end{itemize}
\end{proposition}
\begin{proof}
 \def\UU{\mathcal{U}}
 \def\VV{\mathcal{V}}
 (i) 
 We assume that $U \cap V =\vc{0}$ or, equivalently, $\UU \cap \VV =\emptyset$,
 because otherwise neither $(U,V)$ is a tiling  nor $(\UU,\VV)$ is a factorization.
 Under that assumption, $\vc{0}$ is uniquely represented as the sum
 of an element from $U$ and an element from $V$. 
 Consider a nonzero vector $\vc{y}$ and the $1$-dimensional
 subspace $\pp{y}$ containing it, $\pp{y}=\langle \vc{y}\rangle$.
 It is easy to see the following:
 \begin{itemize}
  \item 
 $\vc{y} = \vc{u}+\vc{0}$, $\vc{u} \in U$, if and only if $\pp y \in \UU$;
  \item 
 $\vc{y} = \vc{0}+\vc{v}$, $\vc{v} \in V$, if and only if $\pp y \in \VV$;
  \item 
 $\vc{y} = \vc{u}+\vc{v}$, $\vc{0} \ne \vc{u} \in U$, $\vc{0} \ne \vc{v} \in V$, 
 if and only if $\pp y \in \plangle \langle \vc{u}\rangle,\langle \vc{v}\rangle\prangle \setminus\{\langle \vc{u}\rangle,\langle \vc{v}\rangle\}$.
%  where $\pp{u} = \langle \vc{u}\rangle$ 
%  and $\pp{v} = \langle \vc{v}\rangle$.
 \end{itemize}
From the three assertions above,
it is now obvious that $(U,V)$ satisfies the definition of a tiling
if and only if $(\UU,\VV)$ satisfies the definition of a factorization.
 
(ii) is straightforward from the obvious
$\plangle \UU \prangle = \PG(\langle U \rangle)$.

(iii) If all vectors from a $1$-dimensional space  
$\pp{x}=\langle \vc{x}\rangle$
are periods of $U$, then for any 
$\vc{y}$ in $U$ it holds 
$\langle \vc{y},\vc{x}\rangle \subset U$, in particular,
$\plangle \langle \vc{y}\rangle,\pp{x}\prangle \subset \UU$.
In this case, by the definition, $\pp{x}$ is a period of~$\UU$.

Inversely, if $\pp{x}=\langle \vc{x}\rangle$ is a period of~$\UU$,
then for any 
$\vc{y}$ in $U\setminus \pp{x}$ ,
the line $\plangle \langle \vc{y}\rangle,\pp{x} \prangle$ is a subset of~$\UU$,
the subspace $\langle \vc{y},\vc{x} \rangle$ is a subset of~$U$,
and any vector from $\pp{x}$ is a period of~$U$.
\end{proof}

The next two propositions are similar to known facts on factorizations of groups.
They emphasize the fundamental value of full-rank and aperiodic factorizations.

\begin{proposition}\label{p:rank}
If $(\mathcal{U},\mathcal{V})$ is a factorization of a projective space
and $\mathcal{U}$ is not full-rank,
then $(\mathcal{U},\plangle \mathcal{U}\prangle \cap \mathcal{V})$
is a factorization of $\plangle \mathcal{U}\prangle$.
\end{proposition}
\begin{proof}
 Let us consider a point $\pp{y}$ in $\plangle \mathcal{U}\prangle$.
 Since $(\mathcal{U},\mathcal{V})$ is a factorization of the whole projective space, exactly one from the following three cases takes place.
 \begin{itemize}
  \item $\pp{y} \in \mathcal{U}$;
  \item $\pp{y} \in \mathcal{V}$; in this case, 
  $\pp{y} \in \plangle \mathcal{U}\prangle \cap \mathcal{V}$;
  \item 
  $\pp{y} \in \plangle \pp{u},\pp{v}\prangle \setminus \{ \pp{u},\pp{v}\}$ 
  for some unique
  $\pp{u}$ from $\mathcal{U}$ and $\pp{v}$ from $\mathcal{V}$;
  in this case, we have $\pp{v} \in \plangle \pp{u},\pp{y}\prangle$,
  and hence $\pp{v} \in \plangle \mathcal{U}\prangle \cap \mathcal{V}$.
 \end{itemize}
 By the definition,
 $(\mathcal{U},\plangle \mathcal{U}\prangle \cap \mathcal{V})$
is a factorization of $\plangle \mathcal{U}\prangle$. 
\end{proof}

For a projective space $\GG$
and a point $\pp{x}$ in it,
by $\GG / \pp{x} $ we denote 
the projective space 
(of the preceding dimension)
whose points
are the lines of  $\GG$ through $\pp{x}$ and lines are induced
by the $2$-dimensional subplanes of $\GG$ through $\pp{x}$
(i.e., $\GG / \pp{x} $
is the ``projection'' of $\GG$
along the lines through~$\pp{x}$).
In particular, 
if 
$\GG = \PG(\FF_q^n)$,
then the points of $\GG / \pp{x} = \PG(\FF_q^n / \pp{x})$
are represented by $2$-dimensional subspaces
of $\FF_q^n$ that include the 
$1$-dimensional subspace~$\pp{x} $,
while the lines of $\GG / \pp{x} $
are represented by $3$-dimensional subspaces
of $\FF_q^n$ including~$\pp{x}$.
For a set $\mathcal C$ of points of~$\GG$,
by $\mathcal C / \pp{x}$ we denote the following set of points of $\GG / \pp{x} $:
$\mathcal C / \pp{x} = \{\plangle \pp{x},\pp{y}\prangle :\ \pp{y} \in \mathcal C,\ \pp{y}\ne \pp{x} \}$.

\begin{proposition}\label{p:period}
If $\GG$ is a projective geometry,
$(\mathcal{U},\mathcal{V})$ is a factorization of $\GG$,
and $\pp{x}$ is a period of $\mathcal{U}$,  
then $(\mathcal{U} / \pp{x},\mathcal{V} / \pp{x})$
is a factorization of $\GG / \pp{x}$.
\end{proposition}
\begin{proof}
 At first, we  prove that 
 \begin{itemize}
  \item[(*)]  \emph{every point of $\GG / \pp{x}$
 belongs to $\mathcal{U} / \pp{x}$,
 or to $\mathcal{V} / \pp{x}$, or to
 $\plangle \ppp{u},\ppp{v} \prangle\setminus\{\ppp{u},\ppp{v} \}$ for some
 $\ppp{u} \in \mathcal{U} / \pp{x}$ and $\ppp{v}\in\mathcal{V} / \pp{x}$.}
 \end{itemize}
 Indeed, consider a point of $\GG / \pp{x}$
 neither in $\mathcal{U} / \pp{x}$
 nor in $\mathcal{V} / \pp{x}$.
 That point, by the definition, is a line through $\pp{x}$ in~$\GG$,
 say $\plangle \pp{x}, \pp{y} \prangle$. 
 Since $(\mathcal{U},\mathcal{V})$ is a factorization of~$\GG$ and 
 $ \pp{y}\not\in \mathcal{U} \cup \mathcal{V} $,
 there are $\pp{u}$ in $\mathcal{U} $ and $\pp{v}$ in $\mathcal{V} $ 
 such that $\pp{y} \in \plangle \pp{u}, \pp{v} \prangle$.
 If $\pp{x} \in \plangle \pp{u}, \pp{v} \prangle$,
 then at least one of~$\pp{u}$, $\pp{v}$ is different from~$\pp{x}$
 and hence $\plangle \pp{x}, \pp{y} \prangle=\plangle \pp{u}, \pp{v} \prangle$  
 belongs to $\mathcal{U} / \pp{x}$ or  $\mathcal{V} / \pp{x}$.
 Otherwise, $\plangle \pp{x}, \pp{y} \prangle$ 
 belongs to the line 
 between  
 the points $\plangle \pp{x}, \pp{u} \prangle$ 
 and $\plangle \pp{x}, \pp{v} \prangle$  in  $\GG / \pp{x}$.
 In any case, (*) holds.
 
 It remains to show that 
 \begin{itemize}
  \item[(**)]  \emph{every  point of $\GG / \pp{x}$
 belongs to only one set in~\rm(*).}
 \end{itemize}
 $\mathcal{U} / \pp{x}$ and $\mathcal{V} / \pp{x}$
 are disjoint because 
 $\mathcal{U} \cap \mathcal{V} = \emptyset$
 and $\pp{x}$ is a period of~$\mathcal{U} $.
 
Now consider the line 
$\plangle \ppp{u},\ppp{v}\prangle$
for some
 $\ppp{u}$ 
%  $\ppp{u}=\plangle \pp{x},\pp{u}\prangle$ 
 in $\mathcal{U} / \pp{x}$ and $\ppp{v}$ 
 in $\mathcal{V} / \pp{x}$.
 The line $\ppp{v}$ of $\GG$ contains a point $\pp{v}$ 
 in~$\mathcal{V}$.
 Moreover, all points of the line $\ppp{u}$ of~$\GG$ 
 belong to $\mathcal{U}$, because $\pp{x}$ is a period of~$\mathcal{U}$.
 In the projective plain 
 $\plangle %\pp{x},
 \ppp{u},\pp{v} \prangle$, 
 a $2$-dimensional subspace of $\GG$,
 every point $\pp{z}$, $\pp{z} \not\in \ppp{u} \cup \ppp{v}$,
 belongs to the line $\plangle\pp{v},\pp{u} \prangle$ between $\pp{v}$ and some point $\pp{u}$
 in $\ppp{u}$. Since $(\mathcal{U},\mathcal{V})$ is a factorization, such point $\pp{z}$ cannot belong to
 $\mathcal{U}$, or $\mathcal{V}$, or any line between
 $\mathcal{U}$ and~$\mathcal{V}$ different from  $\plangle\pp{v},\pp{u} \prangle$ .
 That means that $\plangle\pp{x},\pp{z} \prangle$
 does not belong to $\mathcal{U} / \pp{x}$, or 
 $\mathcal{V} / \pp{x}$, or any line between 
$\mathcal{U} / \pp{x}$ and
 $\mathcal{V} / \pp{x}$
 different from $\plangle \ppp{u},\ppp{v}\prangle$, 
 which completes the proof of~(**).

 By (*) and (**), $(\mathcal{U} / \pp{x},\mathcal{V} / \pp{x})$  satisfies the definition of a factorization.
\end{proof}

\begin{problem}
 Given a projective space,
 characterize all sizes of tiles of full-rank factorizations.
 In particular, which geometries have  no full-rank factorizations
 (in the theory of group factorizations, that property is called 
 the \emph{R\'edei property}~\cite[Ch.\,9]{Szabo:2004:topics})?
 The same question, for aperiodic full-rank factorizations.
\end{problem}

Similar factorizations can be considered for other 
geometries, not only projective spaces.

\begin{problem}
Are there nontrivial ($|\mathcal U|>1$, $|\mathcal V|>1$) examples of factorizations
$(\mathcal U,\mathcal V)$ of the affine spaces, in the sense above
(i.e., when the space is partitioned into the sets
$\mathcal U$, $\mathcal V$, and $\plangle a,b \prangle \backslash \{a,b\}$, 
$a\in \mathcal U$,
$b\in \mathcal V$)?
\end{problem}
For example, in the affine space $\FF_7^3$,
if $|\mathcal U|=5$ and $|\mathcal V|=13$, 
then
$|U|+|V|+\sum_{a\in \mathcal U,b\in \mathcal V}|\plangle a,b \prangle \backslash \{a,b\}|=5+13+5\cdot 13\cdot 5 = 343$,
which is exactly $|\FF_7^3|$.
However, factorizations with such parameters do not exist
(exhaustive search).

Also, factorizations into more than two sets can be considered in a similar manner.
For example, for three sets, such factorization is a partition of the space 
into the sets $\mathcal U$, $\mathcal V$, $\mathcal W$,   
$\plangle a,b \prangle \backslash \{a,b\}$,   
$\plangle a,c \prangle \backslash \{a,c\}$,   
$\plangle b,c \prangle \backslash \{b,c\}$, and
$\plangle a,b,c \prangle 
\backslash \plangle a,b \prangle 
\backslash \plangle a,c \prangle 
\backslash \plangle b,c \prangle $,
$a\in \mathcal U$,
$b\in \mathcal V$,
$c\in \mathcal W$.

\section*{Data Availability Statement}
The dataset generated during the current study 
(examples of full-rank ternary $(13,3^{10},3)_3$ perfect codes) 
is available in the IEEE DataPort repository~\cite{Perfect-related}. 
% \url{https://dx.doi.org/10.21227/w856-4b70}.

% \section*{Compliance with Ethical Standards}
% 
% The author has no competing interests to declare that are relevant to the content of this article.
% For this type of study, informed consent is not required. 
% The work was carried out at the Sobolev Institute of Mathematics at the expense of the Russian Science Foundation, Grant 22-11-00266.

%    \bibliographystyle{abbrv}
%    \bibliography{k}
%   \end{document}

\providecommand\href[2]{#2} \providecommand\url[1]{\href{#1}{#1}}
  \def\DOI#1{{\small {DOI}:
  \href{http://dx.doi.org/#1}{#1}}}\def\DOIURL#1#2{{\small{DOI}:
  \href{http://dx.doi.org/#2}{#1}}}

\end{document}